\documentclass[%
 aip,
 amsmath,amssymb,
preprint,%
]{revtex4-1}

\usepackage{graphicx}
\usepackage{dcolumn}
\usepackage{bm}

\usepackage[utf8]{inputenc}
\usepackage[T1]{fontenc}
\usepackage{mathptmx}
\usepackage{etoolbox}

\usepackage{hyperref}

\makeatletter
\def\@email#1#2{%
 \endgroup
 \patchcmd{\titleblock@produce}
  {\frontmatter@RRAPformat}
  {\frontmatter@RRAPformat{\produce@RRAP{*#1\href{mailto:#2}{#2}}}\frontmatter@RRAPformat}
  {}{}
}%
\makeatother
\begin{document}


\title[Helmholtz Equation on $E_2$]{Reviewing the Helmholtz Equation on Euclidean Plane and Interbasis Expansions}
\author{\framebox{G.S.~Pogosyan}}
 \affiliation{Yerevan State University, Yerevan, Armenia.}
\author{A.~Yakhno}%
 \email{alexander.yakhno@academicos.udg.mx}
\affiliation{%
Departamento de Matematicas, CUCEI, Universidad de Guadalajara, Mexico
}%

\date{\today}

\begin{abstract}
In the present paper we revisit the Helmholtz equation on the Euclidean plane and make some remarks on normalization constants and completeness of wave function sets. The coefficients of interbasis expansions are also reconsidered.
\end{abstract}

\maketitle

%
%
%

\def\vphi{\varphi}
\def\rphi{\vphi\in[0,2\pi)}
\def\vtheta{\vartheta}
\def\SCPM{{\SCPZ\over\SCPN}}
\def\SCPZ{{\xi^2+\eta^2}}
\def\SCPN{{\xi^2\eta^2}}

\def\ralpha{\alpha\in(\i K',\i K'+2K)}
\def\rPalpha{\alpha\in(\i K',\i K'+K)}
\def\rbeta{\beta\in[0,4K')}
\def\rPbeta{\beta\in(0,K')}
\def\rmu{\mu\in(\i K',\i K'+2K)}
\def\rPmu{\mu\in(\i K',\i K'+K)}
\def\reta{\eta\in[0,4K')}
\def\rPeta{\eta\in(0,K')}

\newdimen\theight
\newcommand{\ds}{\displaystyle}
\newcommand{\be}{\begin{equation}}
\newcommand{\ee}{\end{equation}}
\newcommand{\bi}{\begin{itemize}}
\newcommand{\ei}{\end{itemize}}
\newcommand{\x}{{\ensuremath{\times}}}
\newcommand{\bb}[1]{\makebox[16pt]{{\bf#1}}}

\newcommand{\sn}{\mbox{sn}}
\newcommand{\cn}{\mbox{cn}}
\newcommand{\dn}{\mbox{dn}}
\newcommand{\ba}{\begin{array}}
\newcommand{\ea}{\end{array}}
\newcommand{\bea}{\begin{eqnarray}}
\newcommand{\eea}{\end{eqnarray}}
\newcommand{\Res}{\mbox{Res}}
\newcommand{\arcsinh}{\mbox{arsinh}}
\newcommand{\arccosh}{\mbox{arcosh}}
\newcommand{\sign}{\mbox{sign}}
\newcommand{\diag}{\mathop{\mathrm{diag}}}

\def \Column{%
             \vadjust{\setbox0=\hbox{\sevenrm\quad\quad tcol}%
             \theight=\ht0
             \advance\theight by \dp0    \advance\theight by \lineskip
             \kern -\theight \vbox to \theight{\rightline{\rlap{\box0}}%
             \vss}%
             }}%

\catcode`\@=11
\def\qed{\ifhmode\unskip\nobreak\fi\ifmmode\ifinner\else\hskip5\p@\fi\fi
 \hbox{\hskip5\p@\vrule width4\p@ height6\p@ depth1.5\p@\hskip\p@}}
\catcode`@=12 

\def\cents{\hbox{\rm\rlap/c}}
\def\miss{\hbox{\vrule height2pt width 2pt depth0pt}}

\def\vvert{\Vert}                

\def\tcol#1{{\baselineskip=6pt \vcenter{#1}} \Column}

\def\dB{\hbox{{}}}                 
\def\mB#1{\hbox{$#1$}}             
\def\nB#1{\hbox{#1}}               




The Helmholtz equation in the two-dimensional Euclidean space $E_2$ is one of the simplest examples of integrable equations. The classical book of W. Miller \cite{5} discusses the relationship between symmetries, coordinate systems that allow separation of variables and solutions. The process of calculating the coefficients of some interbasis expansions is described in a mathematically rigorous manner.

In this paper we considered plane solutions with a certain parity, their normalization, and also added conditions for the completeness of the sets of the obtained wave functions. This is because on 3D spaces of constant curvature there sometimes arise wave functions with parity that contract to such solutions on $E_2$. It is therefore desirable to collect normalized plane solutions that form a complete basis in order to trace the contraction limit. We also analyzed the properties of the obtained coefficients of interbasis expansions. The presented material is of a reference nature, although it contains a new information. In particular, the coefficients of the expansion of the parabolic basis through the polar one are expressed as polynomials.

The 2d-Helmholtz equation has the form 
\be
\label{EU-01}
\Delta  \Psi + k^2 \Psi = 0,
\qquad
k > 0,
\ee
and allows the separation of variables in four orthogonal coordinate systems: polar, Cartesian, parabolic and elliptic. Each separable coordinate system corresponds to the second-order operator: polar $X_S = L_3^2$, Cartesian $X_C= P_2^2$ and parabolic one $X_P = L_3 P_2 + P_2 L_3$ (here we will not consider the elliptic coordinate system). Operators $L_3 = x\partial_{y} - y\partial_{x}$, $P_1 = \partial_{x}$, $P_2 =  \partial_{y}$ form the basis of $e(2)$ algebra. For $k=0$ we have Laplace equation. In this case, in addition to the four systems mentioned, there are other orthogonal coordinates\cite{5,6} that make it possible to separate the variables.

\section{Solutions}

\subsection{Cartesian coordinates}
For Helmholtz equation in Cartesian coordinates $x, y \in\mathbb{R}$
\be
\label{EU-02}
\frac{\partial^2 \Psi}{\partial x^2} + \frac{\partial^2 \Psi}{\partial y^2} = - k^2\Psi,
\ee
one can separate variables $\Psi(x,y) = X(x)Y(y)$, so
\be
X'' + k_1^2 X = 0, \quad Y'' + k_2^2 Y = 0, \qquad k_1^2 + k_2^2 = k^2, \ k_{1,2}\in\mathbb{R},
\ee
$X_C \Psi = - k_2^2 \Psi$, and formally
\bea
X = C_1 e^{i|k_1|x} + C_2 e^{-i|k_1|x} = A\cos|k_1| x + B\sin|k_1| x, \\
Y = C_3 e^{i|k_2|y} + C_3 e^{-i|k_2|y} = C\cos|k_2| y + D\sin|k_2| y.
\eea
\paragraph{} One can consider the orthonormal complete set of functions
\be
\label{CARTESIAN_SOL}
\Psi_{k_1 k_2}(x,y) = \frac{e^{ik_1 x}}{\sqrt{2\pi}} \frac{e^{ik_2 y}}{\sqrt{2\pi}},
\ee
with normalization condition
\be
\label{NC_E2}
\int\limits_{-\infty}^{\infty} dx  \int\limits_{-\infty}^{\infty} dy \Psi_{k_1 k_2}(x,y) \Psi^\ast_{k_1^\prime k_2^\prime}(x,y)  = \delta(k_1 - k_1^\prime) \delta(k_2 - k_2^\prime),
\ee
and completeness condition
\be
\label{CC_E2}
\int\limits_{-\infty}^{\infty} dk_1  \int\limits_{-\infty}^{\infty} dk_2 \Psi_{k_1 k_2}(x,y) \Psi^\ast_{k_1 k_2}(x^\prime,y^\prime)  = \delta(x - x^\prime) \delta(y - y^\prime),
\ee
here $\delta(t)$ is the Dirac delta-function, the main properties of which can be found, for example, in Refs. \onlinecite{MADELUNG:1957}, \onlinecite{OLVER_F}, \onlinecite{QUANTUM}.
\paragraph{} Considering the parity of wave function with respect to one variable only, we can take the orthonormal complete set of functions in the form (with respect to the change $y \to - y$):
\be
\label{CARTESIAN_P_M}
\Psi^{(+)}_{k_1 k_2}(x,y) = \frac{e^{ik_1 x}}{\sqrt{2\pi}} \frac{\cos |k_2|y}{\sqrt{2\pi}}, \quad 
\Psi^{(-)}_{k_1 k_2}(x,y) = \frac{e^{ik_1 x}}{\sqrt{2\pi}} \frac{\sin |k_2|y}{\sqrt{2\pi}},
\ee
with normalization condition
\be
\label{CARTESIAN_NORMAL_k1k2}
\int\limits_{-\infty}^{\infty} dx  \int\limits_{-\infty}^{\infty} dy \Psi^{(\pm)}_{k_1 k_2}(x,y) \Psi^{(\pm)\ast}_{k_1^\prime k_2^\prime}(x,y)  = \frac{1}{2}\delta(k_1 - k_1^\prime) \delta(|k_2| - |k_2^\prime|),
\ee
and completeness condition
\be
\label{COMLETE_CARTESIAN_k1k2}
\int\limits_{-\infty}^{\infty} dk_1  \int\limits_{-\infty}^{\infty} dk_2 \left[\Psi^{(+)}_{k_1 k_2}(x,y) \Psi^{(+)\ast}_{k_1 k_2}(x^\prime,y^\prime) +  \Psi^{(-)}_{k_1 k_2}(x,y) \Psi^{(-)\ast}_{k_1 k_2}(x^\prime,y^\prime)\right] = \delta(x - x^\prime) \delta(y - y^\prime).
\ee

Let us note, that 
\be
\Psi_{k_1 k_2}(x,y) = \Psi^{(+)}_{k_1 k_2}(x,y) + i\, {\rm sign}(k_2) \Psi^{(-)}_{k_1 k_2}(x,y),
\ee
and relations (\ref{CARTESIAN_NORMAL_k1k2}) and (\ref{COMLETE_CARTESIAN_k1k2}) are in accordance with (\ref{NC_E2}) and (\ref{CC_E2}). Also, one can consider the parity with respect to the change $x \to - x$
\be
\label{CARTESIAN_P_M_x}
\tilde{\Psi}^{(+)}_{k_1 k_2}(x,y) = \cos|k_1|x\, e^{i k_2 y}/2\pi,  \quad \tilde{\Psi}^{(-)}_{k_1 k_2}(x,y) = \sin|k_1|x\, e^{i k_2 y}/2\pi,
\ee
with similar normalization relation.

For parameter $\alpha\in[- \pi, \pi)$ 
\be
\label{K1_K2_K_Alpha}
k_1 = k\cos\alpha, \quad k_2 = k\sin\alpha,
\ee
the right side of (\ref{CARTESIAN_NORMAL_k1k2}) is transformed to
\bea
\frac{\delta(k_1 - k_1^\prime)}{2}\delta\left(\sqrt{k^2 - k_1^2} - \sqrt{k^{\prime2} - k_1^2} \right) = \frac{\delta(k\cos\alpha - k\cos\alpha^\prime)}{2}\frac{|k_2|}{k}\delta(k - k') =  \nonumber\\ 
= \frac{\delta(k - k')}{2k} \delta(|\alpha| - |\alpha'|).
\eea
Therefore functions
\be
\label{CARTESIAN_K_ALPHA_K1_K2}
\Psi^{(\pm)}_{k |\alpha|}(x,y) := \sqrt{k}\, \Psi^{(\pm)}_{k_1 k_2}(x,y), 
\ee
i.e. (taking into account that $|\sin \alpha | = \sin |\alpha|$)
\be
\Psi^{(+)}_{k |\alpha|}(x,y) = \frac{\sqrt{k}}{2\pi} e^{i k x \cos|\alpha|} \cos(k \sin|\alpha|y), \ 
\Psi^{(-)}_{k |\alpha|}(x,y) = \frac{\sqrt{k}}{2\pi} e^{i k x \cos|\alpha|} \sin(k \sin|\alpha|y),
\label{PSI_k_alpha}
\ee
satisfy normalization condition
\be
\int\limits_{- \infty}^{\infty} dx  \int\limits_{- \infty}^{\infty} dy \Psi^{(\pm)}_{k|\alpha|}(x,y) \Psi^{(\pm)\ast}_{k' |\alpha' |}(x,y) = \frac12 \delta(k - k')\delta(|\alpha| - |\alpha'|),
\label{NORMAL_PSI_ALPHA_K}
\ee
while completeness relation (\ref{COMLETE_CARTESIAN_k1k2}) takes the form
\be
\label{COMLETE_CARTESIAN_k_alpha}
\int\limits_{0}^{\infty} dk  \int\limits_{- \pi}^{\pi} d\alpha  \left[\Psi^{(+)}_{k|\alpha|}(x,y) \Psi^{(+)\ast}_{k|\alpha|}(x^\prime,y^\prime) +  \Psi^{(-)}_{k|\alpha|}(x,y) \Psi^{(-)\ast}_{k|\alpha|}(x^\prime,y^\prime)\right] = \delta(x - x^\prime) \delta(y - y^\prime),
\ee
where factor $k$, which follows from Jacobian of transformation (\ref{K1_K2_K_Alpha}), is cancelled because of (\ref{CARTESIAN_K_ALPHA_K1_K2}).
\paragraph{} Considering the parity with respect to both variables, on can take the orthonormal complete set of functions in the form of four different sets:
\bea
\Psi^{(+, +)}_{k_1 k_2}(x,y) =\frac{1}{2\sqrt{\pi}}\cos |k_1|x \cos |k_2|y, \quad  
\Psi^{(+, -)}_{k_1 k_2}(x,y) =\frac{1}{2\sqrt{\pi}}\cos |k_1|x \sin |k_2|y, \\
\Psi^{(-, +)}_{k_1 k_2}(x,y) =\frac{1}{2\sqrt{\pi}}\sin |k_1|x \cos |k_2|y, \quad  
\Psi^{(-, -)}_{k_1 k_2}(x,y) =\frac{1}{2\sqrt{\pi}}\sin |k_1|x \sin |k_2|y,
\eea
with normalization condition
\be
\int\limits_{-\infty}^{\infty} dx  \int\limits_{-\infty}^{\infty} dy \Psi^{(\pm, \pm)}_{k_1 k_2}(x,y) \Psi^{(\pm, \pm)\ast}_{k_1^\prime k_2^\prime}(x,y)  = \frac{1}{4}\delta(|k_1| - |k_1^\prime|) \delta(|k_2| - |k_2|^\prime),
\ee
and completeness condition
\bea
\int\limits_{-\infty}^{\infty} dk_1  \int\limits_{-\infty}^{\infty} dk_2 \left[\Psi^{(+, +)}_{k_1 k_2}(x,y) \Psi^{(+, +)\ast}_{k_1 k_2}(x^\prime,y^\prime) + \Psi^{(+, -)}_{k_1 k_2}(x,y) \Psi^{(+, -)\ast}_{k_1 k_2}(x^\prime,y^\prime) + \right. \nonumber \\
\left.
+ \Psi^{(-, +)}_{k_1 k_2}(x,y) \Psi^{(-, +)\ast}_{k_1 k_2}(x^\prime,y^\prime) + \Psi^{(-, -)}_{k_1 k_2}(x,y) \Psi^{(-, -)\ast}_{k_1 k_2}(x^\prime,y^\prime)
  \right] = \delta(x - x^\prime) \delta(y - y^\prime).
\eea

To verify all integrals we use the following relations for Dirac delta-function:
\bea
\int\limits_{-\infty}^{\infty} e^{it(x - x')} dt=  2\pi \delta(x - x'), \\
\int\limits_{-\infty}^{\infty} \cos at \cos bt \, dt = \pi\left[ \delta(a - b) + \delta(a + b) \right], \\
\int\limits_{-\infty}^{\infty} \sin at \sin bt \, dt = \pi\left[ \delta(a - b) - \delta(a + b) \right].
\eea 

\subsection{Polar coordinates}
In polar coordinates $x = r\cos\vphi$, $y = r\sin\vphi$, $r > 0$, $\vphi\in[0, 2\pi)$ we consider solution in the form
\be
\label{POLAR_WAVE_FUNCTIONS}
\Psi_{km}(r, \vphi) = \sqrt{k} J_{|m|}(kr) \frac{e^{i m\vphi}}{\sqrt{2\pi}}, \qquad  m\in\mathbb{Z},\quad X_S \Psi_{km} = - m^2 \Psi_{km},
\ee
that satisfies normalization condition
\be
\int\limits_{0}^{\infty}r dr \int\limits_{0}^{2\pi} d\vphi \Psi_{km}(r, \vphi)  \Psi^\ast_{k^\prime m^\prime}(r, \vphi) =\delta(k - k^\prime) \delta_{m m^\prime},
\label{NORM_POLAR}
\ee
and completeness relation
\be
\label{COMPLETE_POLAR}
\int\limits_{0}^{\infty} dk \sum\limits_{m = -\infty}^{\infty} \Psi_{km}(r, \vphi)  \Psi^\ast_{k m}(r^\prime, \vphi^\prime) = \frac1r \delta(r - r^\prime) \delta(\vphi - \vphi^\prime).
\ee
To prove the above relations one can use (1.17.21) \cite{OLVER_F} 
\begin{equation}
\label{LEGENDRE-9}
\sum\limits_{m = - \infty}^{\infty} e^{i m (\vphi - \vphi')} = 2 \pi \delta(\vphi - \vphi'),
\end{equation}
and 
\begin{eqnarray}
\label{BESSEL-NORM1}
\int\limits^{\infty}_{0} J_{|m|}(kr)  J_{|m|}(k'r)\, r dr = \frac{1}{k} \delta(k - k')
\end{eqnarray}
with
\begin{eqnarray}
\label{BESSEL-COMPLET1}
\int\limits^{\infty}_{0} J_{|m|}(kr) J_{|m|}(k r')\, k dk = \frac{1}{r} \delta(r - r'),
\end{eqnarray}
which are the well-known orthogonality and completeness relations for Bessel functions on infinite interval (see, for example Ref. \onlinecite{LEON:2014}).

\subsection{Parabolic coordinates}

Let us consider the parabolic coordinate system  
\be
\label{PAR-EU-00}
x = \frac{\xi^2 - \eta^2}{2}, \quad y= \xi \eta, \qquad \xi \ge 0,  \eta \in \mathbb{R},
\ee
where the restriction imposed on $\xi$ guaranties the single-valued transformation between Cartesian and parabolic coordinates. 
Eq. (\ref{EU-02}) takes the form  
\be
\frac{\partial^2 \Psi}{\partial \xi^2} + \frac{\partial^2 \Psi}{\partial \eta^2} + k^2(\xi^2 + \eta^2)\Psi = 0.
\ee
Separation of variables $\Psi = \psi_1(\xi)\psi_2(\eta)$ leads to
\be
\label{PARAB_SEPARATED}
\psi_1^{\prime\prime} + (k^2\xi^2 + 2\beta)\psi_1 = 0, \qquad \psi_2^{\prime\prime} + (k^2\eta^2 - 2\beta)\psi_2 = 0,
\ee
where $\beta\in\mathbb{R}$ is a separation constant and $X_P\Psi = 2\beta \Psi$. The further change  $\tilde{\xi} = \sqrt{2k} \xi$,  $\tilde{\eta} = \sqrt{2k} \eta$ gives
\bea
\label{PAR0-EU-02}
\frac{d^2 \psi_1}{d {\tilde{\xi}}^2} +  \left(\frac{\tilde{\xi}^2}{4} + 
\frac{\beta}{k}\right) \, \psi_1(\tilde{\xi})= 0, 
\qquad
\frac{d^2 \psi_2}{d {\tilde{\eta}}^2} +  \left(\frac{\tilde{\eta}^2}{4} -
\frac{\beta}{k}\right) \, \psi_2(\tilde{\eta})= 0.
\eea
The well known \cite{BE2} solutions of equation 
\bea
\label{PAR1-EU-02}
\frac{d^2 \omega}{d z^2} +  \left(\frac{z^2}{4} - \rho\right) \, \omega= 0, 
\eea
that are real on the real axis can be  chosen as follows 
\bea
\label{PAR1-EU-03}
\omega_1 (z) = e^{-\frac{i z^2}{4}} \,
{_1F_1} \left(\frac14 - \frac{i \rho}{2};  \, \frac12; \,  \frac{i z^2}{2} \right),
\qquad
\omega_2 (z) =  z \, e^{-\frac{i z^2}{4}} \,
{_1F_1} \left(\frac34 - \frac{i \rho}{2};  \, \frac32; \,  \frac{i z^2}{2} \right),
\eea
therefore
\bea
\psi^{(1)}_1(\xi) =  e^{-\frac{i k \xi^2}{2}} \,
{_1F_1} \left(\frac14 + \frac{i \beta}{2k};  \, \frac12; \,  i k \xi^2 \right), \psi^{(2)}_1(\xi) =  \xi e^{-\frac{i k \xi^2}{2}} \,
{_1F_1} \left(\frac34 + \frac{i \beta}{2k};  \, \frac32; \,  i k \xi^2 \right), 
\eea
and
\bea
\psi^{(1)}_2(\eta) =  e^{-\frac{i k \eta^2}{2}} \,
{_1F_1} \left(\frac14 - \frac{i \beta}{2k};  \, \frac12; \,  i k \eta^2 \right), \psi^{(2)}_2(\eta) =  \eta e^{-\frac{i k \eta^2}{2}} \,
{_1F_1} \left(\frac34 - \frac{i \beta}{2k};  \, \frac32; \,  i k \eta^2 \right).
\eea
The right equation of (\ref{PARAB_SEPARATED}) is invariant with respect to the change $\eta \to - \eta$ therefore, one can consider the wave function $\Psi(\xi, \eta)$ with the corresponding parity.  Thus, in parabolic coordinate system, two sets of functions are possible, the even functions $\Psi^{(+)}_{k \beta} (\xi, \eta) = C_{k \beta}^{(+)}\psi_1^{(1)}(\xi) \psi_2^{(1)}(\eta)$ and the odd ones $\Psi^{(-)}_{k \beta} (\xi, \eta) = C_{k \beta}^{(+)}\psi_1^{(2)}(\xi) \psi_2^{(2)}(\eta)$:
\bea
\label{PAR1-EU-04}
\Psi^{(+)}_{k \beta} (\xi, \eta) 
= C_{k \beta}^{(+)}  
e^{- i k\frac{\xi^2+\eta^2}{2}} \,
{_1F_1} \left(\frac14 + \frac{i \beta}{2 k};  \, \frac12; \,  i k \xi^2 \right)
\,
{_1F_1} \left(\frac14 - \frac{i \beta}{2 k};  \, \frac12; \,  i k \eta^2 \right),
\\[3mm]
\label{0-PAR1-EU-04}
\Psi^{(-)}_{k \beta} (\xi, \eta) 
= C_{k \beta}^{(-)}   
 \xi \eta\, e^{- i k\frac{\xi^2+\eta^2}{2}} \,
{_1F_1} \left(\frac34 + \frac{i \beta}{2 k};  \, \frac32; \,  i k \xi^2 \right)
\,
{_1F_1} \left(\frac34 - \frac{i \beta}{2 k};  \, \frac32; \,  i k \eta^2 \right).
\eea
Let us note, that it is possible to take complex conjugation of functions $\psi^{(1,2)}_{1,2}$ to form $\Psi^{(\pm)}_{k \beta} (\xi, \eta)$.


Since functions $\Psi^{(+)}_{k \beta} (\xi, \eta)$ and $\Psi^{(-)}_{k \beta} (\xi, \eta) $ have 
different parity in the variable $\eta$, then 
\bea
\label{PAR1-EU-05}
\int\limits_{0}^{\infty} d \xi  \int\limits_{-\infty}^{\infty} d\eta  \, (\xi^2 + \eta^2)
\Psi^{(\pm)}_{k \beta} (\xi, \eta) \, \Psi^{(\mp)\ast}_{k' \beta'} (\xi, \eta) 
= 0,
\eea
and hence none of the sets is complete separately. 

To calculate normalization constants $C_{k \beta}^{(\pm)}$, satisfying normalization condition
\bea
\label{PAR1-EU-07}
\int\limits_{0}^{\infty} d \xi  \int\limits_{-\infty}^{\infty} d\eta  \, (\xi^2 + \eta^2)
\Psi^{(\pm)}_{k \beta} (\xi, \eta) \, \Psi^{(\pm)\ast}_{k' \beta'} (\xi, \eta) 
= \delta(k-k') \delta(\beta - \beta'),
\eea
one can use the interbasis expansion of parabolic bases in terms of polar bases or Cartesian one.
As a result of the subsections \ref{sub:IE_PARAB_POLAR} we obtain
\bea
\label{PAR1-EU-05}
C_{k \beta}^{(+)}  =  \frac{\left|\Gamma\left(\frac14 + \frac{i \beta}{2 k}\right)\right|^2}{2\sqrt{2}\,\pi^2}, 
\qquad\qquad
C_{k \beta}^{(-)}  =\sqrt{2}\, k \frac{\left|\Gamma\left(\frac34 + \frac{i \beta}{2 k}\right)\right|^2}{\pi^2}.
\eea

\section{Interbasis expansions}

\subsection{Interbasis expansions between Cartesian and polar wave functions}

To find coefficients of expansion
\be
\label{CARTESIAN_POLAR_alpha}
\Psi^{(\pm)}_{k |\alpha|}(x,y) = \sum\limits_{m = -\infty}^{\infty} {\cal S}^{(\pm)\ast}_{k m \alpha}\Psi_{km}(r, \vphi)
\ee
we take into account relations (\ref{CARTESIAN_K_ALPHA_K1_K2}) and orthogonality $\int\limits_{0}^{2\pi} e^{i(m - m')\vphi} d\vphi = 2\pi \delta_{m m'}$.
From  expansion (\ref{CARTESIAN_POLAR_alpha}) we get
\be
\label{S_PLUS}
 {\cal S}^{(+)\ast}_{k m \alpha} J_{|m|}(kr) = \frac{1}{2\pi\sqrt{2\pi}} \int\limits_{0}^{2\pi} e^{i k r \cos\alpha \cos\vphi - im\vphi} \cos(k|\sin\alpha| r\sin\vphi) d\vphi.
\ee

To calculate the integral in (\ref{S_PLUS}) let us consider formula 7.2.4 (27)\cite{BE2}
\be
\label{00-CONT-SPH_EQUI-3}
e^{i z \cos\vphi} = \sum_{m = -\infty}^{\infty} i^m J_m(z) e^{i m \vphi}.
\ee
Then one can obtain 
\bea
\label{00-HORIC-SPHERICAL-01B}
\int\limits_{0}^{2\pi} 
e^{-i |k_1| r \cos\vphi + i k_2 r \sin\vphi}
 e^{- i m\vphi}  d \vphi 
=
\left\{
\begin{array}{c}
\int\limits_{0}^{2\pi}  e^{- i k r \cos(\vphi + \alpha)} \, e^{- i m\vphi} \, d \vphi,\ \cos\alpha > 0, \\[3mm]
\int\limits_{0}^{2\pi}  e^{i k r \cos(\vphi - \alpha)} \, e^{- i m\vphi} \, d \vphi,\ \cos\alpha < 0
\end{array}
\right. 
= 
 \eea
\bea
= \left\{
\begin{array}{c}
  \sum\limits_{n = -\infty}^{\infty} i^{|n|} J_{|n|}(-kr) e^{i n\alpha}\int\limits_{0}^{2\pi} e^{i n \vphi} e^{- i m\vphi}  d \vphi,\ \cos\alpha > 0, \\[3mm]
\sum\limits_{n = -\infty}^{\infty} i^{|n|} J_{|n|}(kr) e^{-i n\alpha}\int\limits_{0}^{2\pi} e^{i n \vphi}  e^{- i m\vphi}  d \vphi,\ \cos\alpha < 0
\end{array}
\right.
= \left\{
\begin{array}{c}
  2\pi (-i)^{|m|} J_{|m|}(kr) e^{i m\alpha},\ \cos\alpha > 0, \\[3mm]
2\pi i^{|m|} J_{|m|}(kr) e^{- i m\alpha},\ \cos\alpha < 0.
\end{array}
\right.
 \nonumber
\eea
In the same way we get 
\bea
\label{00-HORIC-SPHERICAL-01BB}
\int\limits_{0}^{2\pi} \,e^{i |k_1| r \cos\vphi}  e^{i k_2 r \sin\vphi}  e^{- i m\vphi}  d \vphi
= \left\{
\begin{array}{c}
  2\pi i^{|m|} J_{|m|}(kr) e^{- i m\alpha},\ \cos\alpha > 0, \\[3mm]
2\pi (-i)^{|m|} J_{|m|}(kr) e^{i m\alpha},\ \cos\alpha < 0.
\end{array}
\right.
\eea
Therefore
\be
\label{INT_PHI_ALPHA_COS}
\int\limits_{0}^{2\pi} e^{i k r \cos\alpha \cos\vphi - im\vphi} \cos(k|\sin\alpha| r\sin\vphi) d\vphi = 2\pi i^{|m|} J_{|m|}(kr) \cos m\alpha,
\ee
\be
\label{INT_PHI_ALPHA_SIN}
\int\limits_{0}^{2\pi} e^{i k r \cos\alpha \cos\vphi - im\vphi} \sin(k|\sin\alpha| r\sin\vphi) d\vphi = - {\rm sign}(\sin\alpha) 2\pi i^{|m|} J_{|m|}(kr) \sin m\alpha,
\ee
so 
\be
\label{S_P_S_M}
 {\cal S}^{(+)}_{k m \alpha} =  \frac{(- i)^{|m|}}{\sqrt{2\pi}} \cos m \alpha, \quad  {\cal S}^{(-)}_{k m \alpha} = - {\rm sign}(\sin\alpha) \frac{(- i)^{|m|}}{\sqrt{2\pi}} \sin m \alpha.
\ee

It is easy to verify that 
\bea
\label{COMPLETE_Sp_m}
 \int\limits_{-\pi}^{\pi} {\cal S}^{(+)}_{k m \alpha} {\cal S}^{(+)\ast}_{k m' \alpha}  d\alpha = \frac{(- i)^{|m|} i^{|m'|}}{2\pi}  \int\limits_{-\pi}^{\pi} \cos m\alpha \cos m' \alpha\, d\alpha = \frac{\delta_{m m'}}{2},
\eea
\bea
\label{COMPLETE_Sm_m}
 \int\limits_{-\pi}^{\pi} {\cal S}^{(-)}_{k m \alpha} {\cal S}^{(-)\ast}_{k m' \alpha}  d\alpha = \frac{(- i)^{|m|} i^{|m'|}}{2\pi }  \int\limits_{-\pi}^{\pi}  \sin m\alpha \sin m' \alpha\, d\alpha = \frac{\delta_{m m'}}{2}.
\eea

Substitution of expansion (\ref{CARTESIAN_POLAR_alpha}) to the left side of relation (\ref{COMLETE_CARTESIAN_k_alpha}) gives
\bea
\int\limits_{0}^{\infty} dk  \int\limits_{- \pi}^{\pi} d\alpha  \left[\sum\limits_{m, m' = -\infty}^{\infty} {\cal S}^{(+)\ast}_{k m \alpha} {\cal S}^{(+)}_{k m' \alpha}  \Psi_{km}(r, \vphi)  \Psi^\ast_{km'}(r', \vphi')    + \right.  \nonumber \\
\left. + \sum\limits_{m, m' = -\infty}^{\infty} {\cal S}^{(-)\ast}_{k m \alpha} {\cal S}^{(-)}_{k m' \alpha}  \Psi_{km}(r, \vphi)  \Psi^\ast_{km'}(r', \vphi') \right] = \nonumber
\eea
\bea
= \int\limits_{0}^{\infty} dk   \left[\sum\limits_{m, m' = -\infty}^{\infty}\frac{\delta_{m m'}}{2} \Psi_{km}(r, \vphi)  \Psi^\ast_{km'}(r', \vphi')  + \right.  \nonumber \\
\left. + \sum\limits_{m, m' = -\infty}^{\infty}\frac{\delta_{m m'}}{2}  \Psi_{km}(r, \vphi)  \Psi^\ast_{km'}(r', \vphi') \right]  
= \int\limits_{0}^{\infty}  dk  \sum\limits_{m = -\infty}^{\infty}\Psi_{km}(r, \vphi)  \Psi^\ast_{km}(r', \vphi'),
\eea
i.e. the properties (\ref{COMPLETE_Sp_m}), (\ref{COMPLETE_Sm_m}) are in accordance with completeness (\ref{COMPLETE_POLAR}) of polar basis.

Moreover
\bea
\sum\limits_{m = - \infty}^{\infty} {\cal S}^{(\pm)}_{k m \alpha} {\cal S}^{(\pm)\ast}_{k m \alpha'} = \frac{\delta(|\alpha| - |\alpha'|)}{2} = \frac12\left[\delta(\alpha + \alpha') + \delta(\alpha - \alpha')\right],
\eea
and
\be
\sum\limits_{m = - \infty}^{\infty} {\cal S}^{(\pm)}_{k m \alpha} {\cal S}^{(\mp)\ast}_{k m \alpha'} = 0.
\ee

Substitution of expansion
\be
\label{DIRECT_POLAR_CARTESIAN}
\Psi_{km}(r, \vphi) = \int\limits_{-\pi}^{\pi} \left[ {\cal S}^{(+)}_{k m \alpha} \Psi^{(+)}_{k |\alpha|}(x,y)  +  {\cal S}^{(-)}_{k m \alpha} \Psi^{(-)}_{k |\alpha|}(x,y) \right] d\alpha
\ee
to the left side of (\ref{COMPLETE_POLAR}) gives
\bea
\int\limits_{0}^{\infty} dk   \int\limits_{-\pi}^{\pi} d\alpha \int\limits_{-\pi}^{\pi} d\alpha' \sum\limits_{m = -\infty}^{\infty} \left[ {\cal S}^{(+)}_{k m \alpha}  {\cal S}^{(+)\ast}_{k m \alpha'} \Psi^{(+)}_{k |\alpha|}(x,y)\Psi^{(+)\ast}_{k |\alpha'|}(x',y') + \right. \nonumber \\
\left.+  {\cal S}^{(-)}_{k m \alpha}  {\cal S}^{(-)\ast}_{k m \alpha'} \Psi^{(-)}_{k |\alpha|}(x,y)\Psi^{(-)\ast}_{k |\alpha'|}(x',y') \right] =  \nonumber
\eea
\bea
= \int\limits_{0}^{\infty} dk   \int\limits_{-\pi}^{\pi} d\alpha \int\limits_{-\pi}^{\pi} d\alpha'\, \frac{\delta(|\alpha| - |\alpha'|)}{2}  \left[ \Psi^{(+)}_{k |\alpha|}(x,y)\Psi^{(+)\ast}_{k |\alpha'|}(x',y') + \right. \nonumber \\
\left.+ \Psi^{(-)}_{k |\alpha|}(x,y)\Psi^{(-)\ast}_{k |\alpha'|}(x',y') \right] =   \nonumber
\eea
\bea
= \int\limits_{0}^{\infty} dk   \int\limits_{- \pi}^{\pi} d\alpha \left[ \Psi^{(+)}_{k |\alpha|}(x,y)\Psi^{(+)\ast}_{k |\alpha|}(x',y') + \Psi^{(-)}_{k |\alpha|}(x,y)\Psi^{(-)\ast}_{k |\alpha|}(x',y') \right],
\eea
so we obtain the left side of relation (\ref{COMLETE_CARTESIAN_k_alpha}).

Expansion (\ref{DIRECT_POLAR_CARTESIAN}) leads to the well known relation (see Ref. \onlinecite{5}) 
\be
\label{INT_REPRESENT_BESSEL}
2\pi i^{|m|} J_{|m|}(kr) e^{im\vphi} = \int\limits_{-\pi}^{\pi} e^{ikr\cos(\vphi - \alpha)} e^{im\alpha} d\alpha,
\ee
which follows from (\ref{INT_PHI_ALPHA_COS}), (\ref{INT_PHI_ALPHA_SIN}) (or (\ref{00-HORIC-SPHERICAL-01BB})) with change $\vphi \leftrightarrow \alpha$, $\vphi \leftrightarrow \vphi - \pi$ and conjugation.

\subsection{Interbasis expansions between parabolic and polar wave functions}
\label{sub:IE_PARAB_POLAR}

Let us  study the expansion of the parabolic basis (\ref{PAR1-EU-04}) and (\ref{0-PAR1-EU-04}) 
in terms of the polar basis (\ref{POLAR_WAVE_FUNCTIONS})
\be
\label{PARAB-POLAR-01}
\Psi^{ (\pm)}_{k \beta} (\xi, \eta)  =  \sum\limits_{m = -\infty}^{\infty} {\cal W}^{(\pm)}_{k \beta m}
\Psi_{km}(r, \vphi).
\ee
The coordinates are related as follows
\be
\xi^2 = r(1 + \cos\vphi), \qquad \eta^2 = r(1 - \cos\vphi),
\ee
and using orthogonality of $e^{im\vphi}$ on interval $[0, 2\pi)$ we get
\bea
{\cal W}^{(+)}_{k \beta m} \sqrt{2\pi k} J_{|m|}(kr) = C_{k \beta}^{(+)} \int\limits_0^{2\pi} e^{- ikr - im\vphi} \times \nonumber\\
\times {_1F_1} \left(\frac14 + \frac{i \beta}{2 k};  \, \frac12; \,   i k r(1 + \cos\vphi) \right)
\,
{_1F_1} \left(\frac14 - \frac{i \beta}{2 k};  \, \frac12; \,  i k r(1 - \cos\vphi)\right) d\vphi.
\eea
Considering the asymptotics at $r\sim 0$
\be
J_{|m|}(kr) \sim \frac{1}{|m|!} \left(\frac{kr}{2}\right)^{|m|},
\ee
the above relations can be written as follows
\be
\label{INTEGRAL_W}
{\cal W}^{(+)}_{k \beta m} = \frac{C_{k \beta}^{(+)} |m|!}{ \sqrt{2\pi k}} \left(\frac{kr}{2}\right)^{- |m|} \sum\limits_{n, j = 0}^{\infty} (ikr)^{n + j} \frac{\left(\frac14 + \frac{i \beta}{2 k}\right)_n \left(\frac14 - \frac{i \beta}{2 k}\right)_j}{\left(\frac12\right)_n \left(\frac12\right)_j n! j!} I^{(+)}_{nj},
\ee
where we denote
\be
 I^{(+)}_{nj} = \int\limits_0^{2\pi}  (1 + \cos\vphi)^n (1 - \cos\vphi)^j  e^{ - im\vphi}  d\vphi.
\ee
If $n + j > |m|$, then all summands in (\ref{INTEGRAL_W}) goes to zero as $r\sim 0$ due to the presence of the factor $r^{n + j - |m|}$, therefore the only non-zero terms are the terms with $0\le j \le |m| - n$, and $0 \le n \le |m|$.
Let us note, that 
\be
 I^{(+)}_{nj}= 2^{n + j + 1} \int\limits_0^{\pi}  (\cos\vphi)^{2n} (\sin\vphi)^{2j}  e^{ - 2im\vphi}  d\vphi.
\ee
Taking into account  the binomial formula
\bea
\label{EQUIDIST-SPHERIC-17}
(\cos\vphi)^{p} = \frac{1}{2^p} \, \sum_{\ell=0}^{p}\,
\frac{(p)!}{(p-\ell)! (\ell)!} \, e^{i(p-2\ell)\vphi},
\eea
and 1.5.1. (29)\cite{BE1}
\bea
\label{EQUIDIST-SPHERIC-18}
\int_{0}^{\pi} (\sin\vphi)^{\alpha} e^{i \beta \vphi} \, d \vphi
=
\frac{\pi}{2^{\alpha}}\,
\frac{e^{i\frac{\pi}{2} \beta} \, \Gamma(1+\alpha)}
{\Gamma(1+ \frac{\alpha+ \beta}{2}) \Gamma(1+ \frac{\alpha- \beta}{2})},
\qquad \Re(\alpha) > -1,
\eea
we obtain
\bea
 I^{(+)}_{nj}= 
\frac{2^{n + j + 1}}{2^{2n}}  \sum_{\ell=0}^{2n}
\frac{(2n)!}{(2n - \ell)! \ell!}  \int\limits_0^{\pi} 
  (\sin\vphi)^{2j}  e^{i(2n - 2\ell - 2m)\vphi}  d\vphi = \nonumber \\
= \frac{2 \pi (-1)^{n - m}  \Gamma(1 + 2j) }{2^{n + j}}  \sum_{\ell=0}^{2n}
\frac{(2n)!}{(2n - \ell)! \ell!}  \frac{(-1)^\ell }
{\Gamma(1 + j + n - \ell - m) \Gamma(1 +  j - n + \ell + m)}.
\eea
Because of presence of gamma functions, the only nonzero terms in the above sum is with $\ell = 0$, if $m > 0$ and with $\ell = 2n$, if $m < 0$; moreover $n + j = |m|$. Therefore,
\bea
 I^{(+)}_{nj} = \frac{2 \pi (-1)^{n - m}  \Gamma(1 + 2|m| - 2n) }{2^{|m|}} 
\left\{\begin{array}{c}\frac{1}{\Gamma(1 + |m|  - m) \Gamma(1 +  |m| - 2n  + m)},\quad m > 0 \\\frac{1}{\Gamma(1 + |m| - 2n - m) \Gamma(1 +  |m| + m)},\quad m < 0 \end{array}\right. = \nonumber\\
= \frac{2 \pi (-1)^{n - m} }{2^{|m|}}.
\eea
Thus,
\be
\label{W_PLUS}
{\cal W}^{(+)}_{k \beta m} = \frac{C_{k \beta}^{(+)} \sqrt{2\pi} |m|! (- i)^{|m|}}{ \sqrt{ k}}  \sum\limits_{n = 0}^{|m|} \frac{(-1)^n \left(\frac14 + \frac{i \beta}{2 k}\right)_n \left(\frac14 - \frac{i \beta}{2 k}\right)_{|m| - n}}{\left(\frac12\right)_n \left(\frac12\right)_{|m| - n} n! (|m| - n)!} =  \nonumber
\ee
\be
= (- i)^{|m|} C_{k \beta}^{(+)} \sqrt{\frac{2\pi}{k}}  \frac{\left(\frac14 - \frac{i\beta}{2k}\right)_{|m|}}{\left(\frac12\right)_{|m|}}  \,
{_3 F_2}
\left(
\left.\begin{array}{cccc}
\frac12 - |m|  &  - |m| & \frac14 + \frac{i\beta}{2 k}
\\[3mm]
\frac12 & \frac34 + \frac{i\beta}{2k} - |m| & 
\end{array}\right| 1
\right).
\ee
Applying transformation\cite{BAILEY}
\be
{_3 F_2}
\left(
\left.\begin{array}{cccc}
a &  a' & - n
\\[3mm]
c' & 1 - n - c & 
\end{array}\right| 1
\right)
=
\frac{\Gamma(c + a + n)\Gamma(c)}{\Gamma(c + a)\Gamma(c + n)} \,
{_3 F_2}
\left(
\left.\begin{array}{cccc}
a &  c' - a' & - n
\\[3mm]
c' &  c + a & 
\end{array}\right| 1
\right),
\ee
we obtain
\bea
\label{PARAB-POLAR-02}
{\cal W}^{(+)}_{k \beta m} 
&=& \frac{(- i)^{|m|}}{2\sqrt{\pi^3 k}}
\, \left|\Gamma\left(\frac14 + \frac{i \beta}{2 k}\right)\right|^2 \,
{_3 F_2}
\left(
\left.\begin{array}{cccc}
- |m|  &  |m| & \frac14 + \frac{i\beta}{2 k}
\\[3mm]
\frac12 & \frac12,  & 
\end{array}\right| 1
\right),
\eea
if we take $C_{k \beta}^{(+)}$ as in (\ref{PAR1-EU-05}).

For odd functions, we have
\bea
{\cal W}^{(-)}_{k \beta m} \sqrt{2\pi k} J_{|m|}(kr) = C_{k \beta}^{(-)} \int\limits_0^{2\pi} r\sin\vphi\, e^{- ikr - im\vphi} \times \nonumber\\
\times {_1F_1} \left(\frac34 + \frac{i \beta}{2 k};  \, \frac32; \,   i k r(1 + \cos\vphi) \right)
\,
{_1F_1} \left(\frac34 - \frac{i \beta}{2 k};  \, \frac32; \,  i k r(1 - \cos\vphi)\right) d\vphi.
\eea
or
\be
\label{INTEGRAL_W_M}
{\cal W}^{(-)}_{k \beta m} = \frac{C_{k \beta}^{(-)} |m|!}{ \sqrt{2\pi k}} \left(\frac{kr}{2}\right)^{- |m|} \sum\limits_{n, j = 0}^{\infty} r (ikr)^{n + j} \frac{\left(\frac34 + \frac{i \beta}{2 k}\right)_n \left(\frac34 - \frac{i \beta}{2 k}\right)_j}{\left(\frac32\right)_n \left(\frac32\right)_j n! j!} I^{(-)}_{n j},
\ee
\be
I^{(-)}_{n j} = \int\limits_0^{2\pi}  (1 + \cos\vphi)^n (1 - \cos\vphi)^j \sin\vphi\,  e^{ - im\vphi}  d\vphi.
\ee
If $n + j + 1 > |m|$, then all summands in (\ref{INTEGRAL_W_M}) goes to zero as $r\sim 0$ due to the presence of the factor $r^{n + j +1 - |m|}$, therefore we have only terms with $0\le j \le |m| - n - 1$, and $0 \le n \le |m| - 1$. By analogy with $I^{(+)}_{n j}$, we get
\be
I^{(-)}_{n j} = i\pi (-1)^{n + |m|} 2^{1 - |m|},
\ee 
and finally, taking into account (\ref{PAR1-EU-05}),
\bea
\label{PARAB-POLAR-03}
{\cal W}^{(-)}_{k \beta m} 
&=& 2m \frac{(- i)^{|m|}}{\sqrt{\pi^3 k}}
\, \left|\Gamma\left(\frac34 + \frac{i \beta}{2 k}\right)\right|^2 \,
{_3 F_2}
\left(
\left.\begin{array}{cccc}
1 - |m|  &  1 + |m| & \frac34 + \frac{i\beta}{2 k}
\\[3mm]
\frac32 & \frac32  & 
\end{array}\right| 1
\right).
\eea

Thus the interbasis expansions coefficients ${\cal W}^{(\pm)}_{k \beta m}$  are represented 
as polynomials, namely the continuous Hahn polynomials (9.4.1) \cite{KOEKOEK:2010}:
\be
p_n(x;a,b,c,d) = i^n \frac{(a + c)_n (a + d)_n}{n!} \, {_3 F_2}
\left(
\left.\begin{array}{cccc}
- n  &  n + a + b + c + d - 1 & a + ix
\\[3mm]
a + c & a + d  & 
\end{array}\right| 1
\right),
\ee 
so
\be
\label{W_p_H}
{\cal W}^{(+)}_{k \beta m} = \frac{(-1)^{|m|}\, |m|!}{2\sqrt{\pi k}}  \frac{\left|\Gamma\left(\frac14 + \frac{i \beta}{2 k}\right)\right|^2}{ \Gamma^2\left(\frac12 + |m|\right)}\, p_{|m|}\left(\frac{\beta}{2 k};\frac14, \frac14, \frac14, \frac14\right). 
\ee

By analogy, we obtain
\be
\label{W_m_H}
{\cal W}^{(-)}_{k \beta m} = i\, {\rm\, sign }(m)\, \frac{(-1)^{|m|}\, |m|!}{2\sqrt{\pi k}}  \frac{\left|\Gamma\left(\frac34 + \frac{i \beta}{2 k}\right)\right|^2}{ \Gamma^2\left(\frac12 + |m|\right)}\, p_{|m| - 1}\left(\frac{\beta}{2 k};\frac34, \frac34, \frac34, \frac34\right),
\ee
so ${\cal W}^{(\pm)*}_{k \beta m} = \pm {\cal W}^{(\pm)}_{k \beta m}$. 

Using orthogonality relation for continuous Hahn polynomials  (9.4.2) \cite{KOEKOEK:2010}
\bea
\int\limits_{-\infty}^\infty \Gamma(a + ix)\Gamma(b + ix)\Gamma(c - ix)\Gamma(d - ix) p_{n}(x; a, b, c, d) p_{n'}(x; a, b, c, d) dx =  \nonumber \\ 
= 2\pi \frac{\Gamma(n + a + c)\Gamma(n + a + d)\Gamma(n + b + c)\Gamma(n + b + d)}{(2n + a + b + c + d - 1)\Gamma(n + a + b + c + d - 1) n!} \delta_{n n'},
\eea
we get
\be
\label{ORTHO_W}
\int\limits_{-\infty}^\infty{\cal W}^{(+)}_{k \beta m}  {\cal W}^{(+)*}_{k \beta m'} d\beta = \frac{\delta_{|m| |m'|}}{2}, \quad 
\int\limits_{-\infty}^\infty{\cal W}^{(-)}_{k \beta m}  {\cal W}^{(-)*}_{k \beta m'} d\beta = {\rm sign}(mm')\, \frac{\delta_{|m| |m'|}}{2}.
\ee

Using 2.22.2 8\cite{PRUDNIKOV2}, one can obtain integral representation of coefficients ${\cal W}^{(\pm)}_{k \beta m}$:
\bea
{\cal W}^{(+)}_{k \beta m} = \frac{(-i)^{|m|}}{\pi\sqrt{2k}} \int\limits_0^\pi (1 + \cos\vphi)^{-\frac14 - \frac{i\beta}{2k}} (1 - \cos\vphi)^{-\frac14 + \frac{i\beta}{2k}} \cos m\vphi \, d\vphi, \nonumber \\
{\cal W}^{(-)}_{k \beta m} = \frac{(-i)^{|m|}}{\pi\sqrt{2k}} \int\limits_0^\pi (1 + \cos\vphi)^{-\frac14 - \frac{i\beta}{2k}} (1 - \cos\vphi)^{-\frac14 + \frac{i\beta}{2k}} \sin m\vphi\,  d\vphi.
\label{INTEGRAL_REPRESENT_W}
\eea
Then, with the help of (\ref{INTEGRAL_REPRESENT_W}) one can calculate
\be
\label{ORTHOGONAL_SUM_W}
\sum\limits_{m = - \infty}^{\infty} {\cal W}^{(\pm)}_{k \beta m} {\cal W}^{(\pm)*}_{k \beta' m} = \delta(\beta - \beta').
\ee
Indeed, 
\bea
\sum\limits_{m = - \infty}^{\infty} {\cal W}^{(+)}_{k \beta m} {\cal W}^{(+)*}_{k \beta' m} = \frac{1}{2k\pi^2} \int\limits_0^\pi \frac{d\vphi }{\sqrt{|\sin\vphi|}} \left(\frac{1 - \cos\vphi}{1 + \cos\vphi}\right)^{\frac{i\beta}{2k}}  \int\limits_0^\pi \frac{d\vphi' }{\sqrt{|\sin\vphi'|}} \left(\frac{1 + \cos\vphi'}{1 - \cos\vphi'}\right)^{\frac{i\beta'}{2k}} \times \nonumber\\
\times \sum\limits_{m = - \infty}^{\infty} \cos m\vphi \, \cos m\vphi'.
\eea
Taking into account, that
\be
 \sum\limits_{m = - \infty}^{\infty} \cos m\vphi \, \cos m\vphi' = \pi \delta(\vphi - \vphi'), \quad \vphi,\vphi' \in (0,\pi),
 \label{COS_DELTA}
\ee
we obtain
\be
\sum\limits_{m = - \infty}^{\infty} {\cal W}^{(+)}_{k \beta m} {\cal W}^{(+)*}_{k \beta' m} =  \frac{1}{2k\pi} \int\limits_0^\pi \frac{d\vphi }{\sin\vphi} \left(\frac{1 - \cos\vphi}{1 + \cos\vphi}\right)^{\frac{i(\beta - \beta')}{2k}}.
\ee
Change of variables $\cos\vphi = \tanh\tau$ leads to the integral
\be
\frac{1}{2k\pi} \int\limits_{-\infty}^{\infty} e^{-i(\beta - \beta')\tau/k} d\tau = \delta(\beta - \beta').
\label{INT_beta}
\ee

For coefficients ${\cal W}^{(-)}_{k \beta m}$ we make the same steps, with the only difference, that we use
\be
\sum\limits_{m = - \infty}^{\infty} \sin m\vphi \, \sin m\vphi' = \pi \delta(\vphi - \vphi').
\label{SIN_SUM}
\ee 
The presence of ${\rm \sign}(m)$ in (\ref{W_m_H}) permits to conclude that
\be
\sum\limits_{m = - \infty}^{\infty} {\cal W}^{(\pm)}_{k \beta m} {\cal W}^{(\mp)*}_{k \beta' m} = 0.
\ee

Let us note that the selection of constants $C^{(\pm)}_{k\beta}$ (\ref{PAR1-EU-05}) and properties (\ref{ORTHOGONAL_SUM_W}) permit to demonstrate orthogonality relation (\ref{PAR1-EU-07}) for parabolic wave functions using decomposition (\ref{PARAB-POLAR-01}) and orthogonality condition (\ref{NORM_POLAR}) for polar basis.

Properties of coefficients ${\cal W}^{(\pm)}_{k \beta m}$ give the inverse decomposition
\be
\label{INTER_POLAR_PARABOLIC}
\Psi_{km}(r, \vphi) = \int\limits_{-\infty}^{\infty}\left[ {\cal W}^{(+)*}_{k \beta m}\Psi^{(+)}_{k \beta} (\xi, \eta)  +
 {\cal W}^{(-)*}_{k \beta m}\Psi^{(-)}_{k \beta} (\xi, \eta)\right] d\beta,
\ee
and also the completeness relation for parabolic basis
\bea
\int\limits_{0}^{\infty} dk \int\limits_{-\infty}^{\infty} d\beta \left[\Psi^{(+)}_{k \beta} (\xi, \eta) \Psi^{(+)*}_{k \beta} (\xi', \eta')  + \Psi^{(-)}_{k \beta} (\xi, \eta) \Psi^{(-)*}_{k \beta} (\xi', \eta') \right] = \nonumber
\eea
\bea
= \int\limits_{0}^{\infty} dk \int\limits_{-\infty}^{\infty} d\beta \left[\sum\limits_{m,m' = 0}^{\infty}  {\cal W}^{(+)}_{k \beta m}  {\cal W}^{(+)*}_{k \beta' m} \Psi_{k m} (r, \vphi) \Psi^{*}_{k m'} (r', \vphi') + \right. \nonumber \\
\left. + \sum\limits_{m,m' = 0}^{\infty}  {\cal W}^{(-)}_{k \beta m}  {\cal W}^{(-)*}_{k \beta' m} \Psi_{k m} (r, \vphi) \Psi^{*}_{k m'} (r', \vphi')
\right] =  \nonumber
\eea
\bea
=   \int\limits_{0}^{\infty} dk \left[ \sum\limits_{m,m' = 0}^{\infty} \frac{\delta_{mm'} + \delta_{m,-m'}}{2}  \Psi_{k m} (r, \vphi) \Psi^{*}_{k m'} (r', \vphi') + \sum\limits_{m,m' = 0}^{\infty} \frac{\delta_{mm'} - \delta_{m,-m'}}{2}  \Psi_{k m} (r, \vphi) \Psi^{*}_{k m'} (r', \vphi')\right] =  \nonumber
\eea
\bea
=  \int\limits_{0}^{\infty} dk \sum\limits_{m = - \infty}^{\infty}  \Psi_{k m} (r, \vphi) \Psi^{*}_{k m} (r', \vphi') = \frac{\delta(\xi - \xi')\delta(|\eta| - |\eta'|) }{\xi^2 + \eta^2},
\label{COMPLETE_PARABOLIC}
\eea
where we use (\ref{COMPLETE_POLAR}) and that
\bea
\frac1r \delta(r - r^\prime) \delta(\vphi - \vphi^\prime) = \frac{2 \delta(\vphi - \vphi') }{\xi^2 + \eta^2}\delta\left(\frac{\xi^2}{1 + \cos\vphi} - \frac{{\xi'}^2}{1 + \cos\vphi} \right) =  \nonumber\\
= \frac{1 + \cos\vphi}{\xi(\xi^2 + \eta^2)} \delta(\xi - \xi') \delta\left(\arccos\frac{\xi^2 - \eta^2}{\xi^2 + \eta^2} - \arccos\frac{\xi^2 - {\eta'}^2}{\xi^2 + {\eta'}^2}\right) = \nonumber\\
= \frac{1 + \cos\vphi}{2\xi^2} \delta(\xi - \xi') \left(\delta(\eta - \eta') + \delta(\eta + \eta') \right) = \nonumber\\
= \frac{\delta(\xi - \xi')\delta(|\eta| - |\eta'|) }{\xi^2 + \eta^2}.
\eea

\subsection{Interbasis expansions between parabolic and Cartesian wave functions}
\label{subsec:Parabolic_Cartesian}

To obtain interbasis expansion of parabolic basis in terms of Cartesian, let us use expansions (\ref{PARAB-POLAR-01}) and (\ref{DIRECT_POLAR_CARTESIAN}). Thus,
\bea
\Psi^{ (\pm)}_{k \beta} (\xi, \eta)  =  \sum\limits_{m = -\infty}^{\infty} {\cal W}^{(\pm)}_{k \beta m}
\Psi_{km}(r, \vphi) = \nonumber\\
=   \int\limits_{-\pi}^{\pi} d\alpha \left[ \Psi^{(+)}_{k |\alpha|}(x,y) \sum\limits_{m = -\infty}^{\infty} {\cal W}^{(\pm)}_{k \beta m}  {\cal S}^{(+)}_{k m \alpha}   +  \Psi^{(-)}_{k |\alpha|}(x,y) \sum\limits_{m = -\infty}^{\infty} {\cal W}^{(\pm)}_{k \beta m}  {\cal S}^{(-)}_{k m \alpha}  \right] = \nonumber \\
= \int\limits_{-\pi}^{\pi} {\cal Z}^{(\pm)}_{k\beta \alpha} \Psi^{(\pm)}_{k |\alpha|}(x,y)  d\alpha,
\label{PARABOLIC_CARTESIAN}
\eea
where we denote
\be
 {\cal Z}^{(\pm)}_{k\beta \alpha} = \sum\limits_{m = -\infty}^{\infty} {\cal W}^{(\pm)}_{k \beta m}  {\cal S}^{(\pm)}_{k m \alpha},
 \label{SUM_W_S}
\ee
and note that
\be
\sum\limits_{m = -\infty}^{\infty} {\cal W}^{(\pm)}_{k \beta m}  {\cal S}^{(\mp)}_{k m \alpha} = 0
\ee
due to presence of  $\sin m\alpha$ in ${\cal S}^{(-)}_{k m \alpha}$ (\ref{S_P_S_M}), and ${\rm\, sign }\, (m)$ in ${\cal W}^{(-)}_{k \beta m} $ (\ref{W_m_H}).

It is easy to see, that functions $\Psi^{(\pm)}_{k |\alpha|}(x,y)$ (\ref{PSI_k_alpha}) and coefficients $ {\cal S}^{(\pm)}_{k m \alpha}$  (\ref{S_P_S_M}) are even functions with respect to change $\alpha \to - \alpha$, therefore coefficients $ {\cal Z}^{(\pm)}_{k\beta \alpha}$ are even functions too, and expansion (\ref{PARABOLIC_CARTESIAN}) takes the form
\be
\Psi^{ (\pm)}_{k \beta} (\xi, \eta)  = 2 \int\limits_{0}^{\pi} {\cal Z}^{(\pm)}_{k\beta |\alpha|} \Psi^{(\pm)}_{k |\alpha|}(x,y)  d\alpha.
\label{PARABOLIC_CARTESIAN2}
\ee

Using integral representation (\ref{INTEGRAL_REPRESENT_W}), from (\ref{SUM_W_S}) and (\ref{S_P_S_M}) we get 
\bea
 {\cal Z}^{(+)}_{k\beta |\alpha|} = \frac{1}{2\pi \sqrt{\pi k}}  \int\limits_0^\pi (1 + \cos\vphi)^{-\frac14 - \frac{i\beta}{2k}} (1 - \cos\vphi)^{-\frac14 + \frac{i\beta}{2k}} \, d\vphi \sum\limits_{m = -\infty}^{\infty} (-1)^{|m|} \cos m \alpha \cos m\vphi. 
\eea
Using that $(-1)^{|m|} \cos m \alpha = \cos(m (\pi -  |\alpha|))$ and (\ref{COS_DELTA}), we obtain
\be
\sum\limits_{m = -\infty}^{\infty} (-1)^{|m|} \cos m \alpha \cos m\vphi = \pi \delta(\vphi - (\pi - |\alpha|)).
\ee
If $t\in(a,b)$, we have (see Appendix II, (6)\cite{QUANTUM})
\be
\int\limits_a^b f(x) \delta(x - t) dx = f(t) 
\ee
and the above integral is equal to zero, if $t\notin(a,b)$. Therefore, for all $|\alpha| \in (0,\pi)$ we obtain 
\bea
 {\cal Z}^{(+)}_{k\beta |\alpha|} = \frac{ (1 + \cos(\pi - |\alpha|))^{-\frac14 - \frac{i\beta}{2k}} (1 - \cos(\pi - |\alpha|))^{-\frac14 + \frac{i\beta}{2k}} }{2 \sqrt{\pi k}} = \frac{1}{2\sqrt{\pi k \sin|\alpha|}} \left(\cot\frac{|\alpha|}{2}\right)^{\frac{i\beta}{k}}.
\eea

By analogy, for odd coefficients with the help of (\ref{SIN_SUM}), we obtain the same expression
\bea
 {\cal Z}^{(-)}_{k\beta |\alpha|} =   {\cal Z}^{(+)}_{k\beta |\alpha|}  = \frac{1}{2\sqrt{\pi k \sin|\alpha|}} \left(\cot\frac{|\alpha|}{2}\right)^{\frac{i\beta}{k}}, \qquad |\alpha| \in (0, \pi).
\eea

Substitution of  (\ref{W_p_H}), (\ref{W_m_H}) to (\ref{SUM_W_S}) formally gives relation between two particular Hahn polynomials
\bea
\left|\Gamma\left(\frac14 + \frac{i \beta}{2 k}\right)\right|^2 \sum\limits_{m = -\infty}^{\infty} \frac{i^{|m|}  |m|!}{ \Gamma^2\left(\frac12 + |m|\right)}  p_{|m|}\left(\frac{\beta}{2 k};\frac14, \frac14, \frac14, \frac14\right) \cos |m|\alpha = \nonumber \\
 = - i \left|\Gamma\left(\frac34 + \frac{i \beta}{2 k}\right)\right|^2 \sum\limits_{m = -\infty}^{\infty} \frac{i^{|m|}  |m|!}{ \Gamma^2\left(\frac12 + |m|\right)}  p_{|m| - 1}\left(\frac{\beta}{2 k};\frac34, \frac34, \frac34, \frac34\right) \sin |m|\alpha = \nonumber \\
 = \frac{\sqrt{2\pi}}{\sqrt{\sin|\alpha|}} \left(\cot\frac{|\alpha|}{2}\right)^{\frac{i\beta}{k}}, \qquad |\alpha|\in(0, \pi).
\eea


Let us deduce some properties of coefficients ${\cal Z}^{(\pm)}_{k\beta |\alpha|}$. Firstly, one can see that
\be
\int\limits_{-\infty}^{\infty} {\cal Z}^{(\pm)}_{k\beta |\alpha|} {\cal Z}^{(\pm)*}_{k\beta |\alpha'|} d\beta = \frac{1}{4\pi k \sqrt{|\sin\alpha \sin\alpha'|}} \int\limits_{-\infty}^{\infty} \left(\frac{\cot\frac{|\alpha|}{2}}{\cot\frac{|\alpha'|}{2}}\right)^{\frac{i\beta}{k}} d\beta.
\ee
Taking into account, that 
\bea
\int\limits_{-\infty}^{\infty} \left(\frac{\cot\frac{|\alpha|}{2}}{\cot\frac{|\alpha'|}{2}}\right)^{\frac{i\beta}{k}} d\beta = \left. - ik \frac{\left(\frac{\cot\frac{|\alpha|}{2}}{\cot\frac{|\alpha'|}{2}}\right)^{\frac{i\beta}{k}}}{\ln \cot\frac{|\alpha|}{2} - \ln \cot\frac{|\alpha'|}{2}} \right|^{\infty}_{-\infty} = 2\pi k \delta\left( \ln \cot\frac{|\alpha|}{2} - \ln \cot\frac{|\alpha'|}{2} \right) = \nonumber \\
=2\pi k |\sin\alpha|\, \delta(|\alpha| - |\alpha'|),
\eea
we obtain
\be
\int\limits_{-\infty}^{\infty} {\cal Z}^{(\pm)}_{k\beta |\alpha|} {\cal Z}^{(\pm)*}_{k\beta |\alpha'|} d\beta = \frac{ \delta(|\alpha| - |\alpha'|)}{2}.
\label{Z_beta}
\ee
Moreover
\be
\int\limits_{-\pi}^{\pi} {\cal Z}^{(\pm)}_{k\beta |\alpha|} {\cal Z}^{(\pm)*}_{k\beta' |\alpha|} d\alpha =  \frac{1}{2\pi k} \int\limits_{0}^{\pi} \frac{d\alpha}{\sin\alpha}  \left(\cot\frac{\alpha}{2}\right)^{\frac{i(\beta - \beta')}{k}} = \delta(\beta - \beta'),
\label{Z_delta_beta}
\ee
where the integral in the above expression is similar to (\ref{INT_beta}).

Relation (\ref{Z_beta}) leads to inverse expansion
\be
\Psi^{(\pm)}_{k |\alpha|}(x,y) =  \int\limits_{-\infty}^{\infty} {\cal Z}^{(\pm)\ast}_{k\beta |\alpha|} \Psi^{ (\pm)}_{k \beta} (\xi, \eta)  d\beta.
\label{INV_CARTESIAN_PARABOLIC}
\ee
Using this expansion and property (\ref{Z_delta_beta}) in left side of (\ref{COMLETE_CARTESIAN_k_alpha}), we obtain
\bea
\int\limits_{0}^{\infty} dk  \int\limits_{- \pi}^{\pi} d\alpha  \left[\int\limits_{-\infty}^{\infty} {\cal Z}^{(+)\ast}_{k\beta |\alpha|} \Psi^{ (+)}_{k \beta} (\xi, \eta)  d\beta  \int\limits_{-\infty}^{\infty} {\cal Z}^{(+)}_{k\beta' |\alpha|} \Psi^{(+)\ast}_{k \beta'} (\xi', \eta')  d\beta' + \right. \nonumber \\ 
\left. + \int\limits_{-\infty}^{\infty} {\cal Z}^{(-)\ast}_{k\beta |\alpha|} \Psi^{ (-)}_{k \beta} (\xi, \eta)  d\beta  \int\limits_{-\infty}^{\infty} {\cal Z}^{(-)}_{k\beta' |\alpha|} \Psi^{(-)\ast}_{k \beta'} (\xi', \eta')  d\beta'\right] = \nonumber
\eea
\bea
= \int\limits_{0}^{\infty} dk \left[
\int\limits_{-\infty}^{\infty} d\beta  \Psi^{ (+)}_{k \beta} (\xi, \eta) \int\limits_{-\infty}^{\infty} d\beta'  \Psi^{ (+)*}_{k \beta'} (\xi', \eta') \delta(\beta - \beta') + \right. \nonumber \\
+ \left. \int\limits_{-\infty}^{\infty} d\beta  \Psi^{ (-)}_{k \beta} (\xi, \eta) \int\limits_{-\infty}^{\infty} d\beta'  \Psi^{ (-)*}_{k \beta'} (\xi', \eta') \delta(\beta - \beta')
 \right] = \frac{\delta(\xi - \xi')\delta(|\eta| - |\eta'|) }{\xi^2 + \eta^2},
\eea
where we apply equality (\ref{COMPLETE_PARABOLIC}). Taking into account, that
\bea
\frac{\delta(\xi - \xi')\delta(|\eta| - |\eta'|) }{\xi^2 + \eta^2} = \frac{\delta(\xi - \xi')}{2\sqrt{x^2 + y^2}}\delta\left(\frac{|y|}{\xi} -  \frac{|y'|}{\xi}\right) = \frac{\xi \delta(|y| - |y'|)}{2\sqrt{x^2 + y^2}} \times \nonumber \\
\times \delta\left(\sqrt{\sqrt{x^2 + y^2} - x} - \sqrt{\sqrt{{x'}^2 + y^2} - x'}\right) 
= \delta(x - x^\prime) \delta(|y| - |y^\prime|),
\eea
we get the right side of completeness relation for Cartesian basis $\Psi^{(\pm)}_{k |\alpha|}(x,y)$ if $yy' \ge 0$.

Substitution of expansion (\ref{PARABOLIC_CARTESIAN}) to the left side of (\ref{PAR1-EU-07}) gives
\bea
\frac{1}{\pi\sqrt{k k'}} \int\limits_{0}^{\pi} \frac{d\alpha}{\sqrt{\sin\alpha}} \left(\cot\frac{\alpha}{2}\right)^{\frac{i\beta}{k}}  \int\limits_{0}^{\pi} \frac{d\alpha'}{\sqrt{\sin\alpha'}} \left(\cot\frac{\alpha'}{2}\right)^{\frac{- i\beta'}{k'}} 
\int\limits_{- \infty}^{\infty} dx  \int\limits_{- \infty}^{\infty} dy \Psi^{(\pm)}_{k|\alpha|}(x,y) \Psi^{(\pm)\ast}_{k' |\alpha' |}(x,y) = \nonumber
\eea
\bea
= \frac{\delta(k - k')}{2 \pi k} \int\limits_{0}^{\pi} \frac{d\alpha}{\sin\alpha} \left(\cot\frac{\alpha}{2}\right)^{\frac{i(\beta - \beta')}{k}} = \delta(k - k') \delta(\beta - \beta'),
\eea
where we use relation (\ref{INT_beta}) and (\ref{NORMAL_PSI_ALPHA_K}). Thus, we again come to the right side of the orthogonality relation (\ref{PAR1-EU-07}).

The parabolic basis ${\Psi}^{(3)}_{\mu\pm}(\xi, \eta)$ given in Ref. \onlinecite{5} is connected to our according to 
\[
{\Psi}^{(3)}_{\mu\pm}(\xi, \eta) = \pi\sqrt{2}\left(\Psi^{(+)}_{k \beta} (\xi, \eta)
\pm i \Psi^{(-)}_{k \beta} (\xi, \eta) \right),
\]
using relation of parabolic cylinder function $D_\nu(z)$ with confluent hypergeometric function $_1F_1$ 8.1.1\cite{MAGNUS}.  Thus, decomposition (\ref{INV_CARTESIAN_PARABOLIC}) is in accordance with (3.49) (Ch. 1 \cite{5}). Note,  that formulas (\ref{W_p_H}) and (\ref{W_m_H}) are more convenient than those given in book \onlinecite{5} since they are expressed in terms of polynomials.


\nocite{*}
\bibliography{bib_Georgy}

\providecommand{\noopsort}[1]{}\providecommand{\singleletter}[1]{#1}%
\begin{thebibliography}{12}%
\makeatletter
\providecommand \@ifxundefined [1]{%
 \@ifx{#1\undefined}
}%
\providecommand \@ifnum [1]{%
 \ifnum #1\expandafter \@firstoftwo
 \else \expandafter \@secondoftwo
 \fi
}%
\providecommand \@ifx [1]{%
 \ifx #1\expandafter \@firstoftwo
 \else \expandafter \@secondoftwo
 \fi
}%
\providecommand \natexlab [1]{#1}%
\providecommand \enquote  [1]{``#1''}%
\providecommand \bibnamefont  [1]{#1}%
\providecommand \bibfnamefont [1]{#1}%
\providecommand \citenamefont [1]{#1}%
\providecommand \href@noop [0]{\@secondoftwo}%
\providecommand \href [0]{\begingroup \@sanitize@url \@href}%
\providecommand \@href[1]{\@@startlink{#1}\@@href}%
\providecommand \@@href[1]{\endgroup#1\@@endlink}%
\providecommand \@sanitize@url [0]{\catcode `\\12\catcode `\$12\catcode
  `\&12\catcode `\#12\catcode `\^12\catcode `\_12\catcode `\%12\relax}%
\providecommand \@@startlink[1]{}%
\providecommand \@@endlink[0]{}%
\providecommand \url  [0]{\begingroup\@sanitize@url \@url }%
\providecommand \@url [1]{\endgroup\@href {#1}{\urlprefix }}%
\providecommand \urlprefix  [0]{URL }%
\providecommand \Eprint [0]{\href }%
\providecommand \doibase [0]{http://dx.doi.org/}%
\providecommand \selectlanguage [0]{\@gobble}%
\providecommand \bibinfo  [0]{\@secondoftwo}%
\providecommand \bibfield  [0]{\@secondoftwo}%
\providecommand \translation [1]{[#1]}%
\providecommand \BibitemOpen [0]{}%
\providecommand \bibitemStop [0]{}%
\providecommand \bibitemNoStop [0]{.\EOS\space}%
\providecommand \EOS [0]{\spacefactor3000\relax}%
\providecommand \BibitemShut  [1]{\csname bibitem#1\endcsname}%
\let\auto@bib@innerbib\@empty
\bibitem [{\citenamefont {Miller~Jr}(1977)}]{5}%
  \BibitemOpen
  \bibfield  {author} {\bibinfo {author} {\bibfnamefont {W.}~\bibnamefont
  {Miller~Jr}},\ }\href@noop {} {\emph {\bibinfo {title} {Symmetry and
  Separation of Variables}}}\ (\bibinfo  {publisher} {Addison Wesley Publishing
  Company},\ \bibinfo {year} {1977})\BibitemShut {NoStop}%
\bibitem [{\citenamefont {Kalnins}(1986)}]{6}%
  \BibitemOpen
  \bibfield  {author} {\bibinfo {author} {\bibfnamefont {E.}~\bibnamefont
  {Kalnins}},\ }\href@noop {} {\emph {\bibinfo {title} {Separation of Variables
  for Riemannian Spaces of Constant Curvature}}}\ (\bibinfo  {publisher}
  {Longman},\ \bibinfo {year} {1986})\BibitemShut {NoStop}%
\bibitem [{\citenamefont {Madelung}(1957)}]{MADELUNG:1957}%
  \BibitemOpen
  \bibfield  {author} {\bibinfo {author} {\bibfnamefont {E.}~\bibnamefont
  {Madelung}},\ }\enquote {\bibinfo {title} {Die mathematischen hilfsmittel des
  physikers},}\ \ (\bibinfo  {publisher} {Springer Berlin Heidelberg},\
  \bibinfo {year} {1957})\ Chap.\ \bibinfo {chapter} {Zahlen, Funktionen und
  Operatoren}\BibitemShut {NoStop}%
\bibitem [{\citenamefont {Olver}\ \emph {et~al.}(2010)\citenamefont {Olver},
  \citenamefont {Lozier}, \citenamefont {Boisvert},\ and\ \citenamefont
  {Clark}}]{OLVER_F}%
  \BibitemOpen
  \bibfield  {author} {\bibinfo {author} {\bibfnamefont {F.}~\bibnamefont
  {Olver}}, \bibinfo {author} {\bibfnamefont {D.}~\bibnamefont {Lozier}},
  \bibinfo {author} {\bibfnamefont {R.}~\bibnamefont {Boisvert}}, \ and\
  \bibinfo {author} {\bibfnamefont {C.}~\bibnamefont {Clark}},\ }\href@noop {}
  {\emph {\bibinfo {title} {NIST Handbook of Mathematical Functions}}}\
  (\bibinfo  {publisher} {Cambridge University Press},\ \bibinfo {year}
  {2010})\BibitemShut {NoStop}%
\bibitem [{\citenamefont {Cohen-Tannoudji}, \citenamefont {Diu},\ and\
  \citenamefont {Lalo\"e}(2020)}]{QUANTUM}%
  \BibitemOpen
  \bibfield  {author} {\bibinfo {author} {\bibfnamefont {C.}~\bibnamefont
  {Cohen-Tannoudji}}, \bibinfo {author} {\bibfnamefont {B.}~\bibnamefont
  {Diu}}, \ and\ \bibinfo {author} {\bibfnamefont {F.}~\bibnamefont
  {Lalo\"e}},\ }\href@noop {} {\emph {\bibinfo {title} {Quantum Mechanics}}},\
  \bibinfo {edition} {2nd}\ ed.,\ Vol.~\bibinfo {volume} {II}\ (\bibinfo
  {publisher} {Wiley-VCH},\ \bibinfo {year} {2020})\BibitemShut {NoStop}%
\bibitem [{\citenamefont {Ponce~de Leon}(2015)}]{LEON:2014}%
  \BibitemOpen
  \bibfield  {author} {\bibinfo {author} {\bibfnamefont {J.}~\bibnamefont
  {Ponce~de Leon}},\ }\bibfield  {title} {\enquote {\bibinfo {title}
  {Revisiting the orthogonality of bessel functions of the first kind on an
  infinite interval},}\ }\href@noop {} {\bibfield  {journal} {\bibinfo
  {journal} {Eur. J. Phys.}\ ,\ \bibinfo {pages} {36 015016}} (\bibinfo {year}
  {2015})}\BibitemShut {NoStop}%
\bibitem [{\citenamefont {Bateman}\ and\ \citenamefont
  {Erd{\'e}lyi}(1952)}]{BE2}%
  \BibitemOpen
  \bibfield  {author} {\bibinfo {author} {\bibfnamefont {H.}~\bibnamefont
  {Bateman}}\ and\ \bibinfo {author} {\bibfnamefont {A.}~\bibnamefont
  {Erd{\'e}lyi}},\ }\href@noop {} {\emph {\bibinfo {title} {Higher
  Transcendental Functions}}},\ Vol.~\bibinfo {volume} {2}\ (\bibinfo
  {publisher} {Mc Graw-Hill},\ \bibinfo {address} {New York},\ \bibinfo {year}
  {1952})\BibitemShut {NoStop}%
\bibitem [{\citenamefont {Bateman}\ and\ \citenamefont
  {Erd{\'e}lyi}(1953)}]{BE1}%
  \BibitemOpen
  \bibfield  {author} {\bibinfo {author} {\bibfnamefont {H.}~\bibnamefont
  {Bateman}}\ and\ \bibinfo {author} {\bibfnamefont {A.}~\bibnamefont
  {Erd{\'e}lyi}},\ }\href@noop {} {\emph {\bibinfo {title} {Higher
  Transcendental Functions}}},\ Vol.~\bibinfo {volume} {1}\ (\bibinfo
  {publisher} {Mc Graw-Hill},\ \bibinfo {address} {New York},\ \bibinfo {year}
  {1953})\BibitemShut {NoStop}%
\bibitem [{\citenamefont {Bailey}(1964)}]{BAILEY}%
  \BibitemOpen
  \bibfield  {author} {\bibinfo {author} {\bibfnamefont {W.}~\bibnamefont
  {Bailey}},\ }\href@noop {} {\emph {\bibinfo {title} {Generalized
  Hypergeometric Series}}}\ (\bibinfo  {publisher} {Stechert-Hafner Service
  Agency},\ \bibinfo {address} {New York},\ \bibinfo {year} {1964})\BibitemShut
  {NoStop}%
\bibitem [{\citenamefont {Koekoek}, \citenamefont {Lesky},\ and\ \citenamefont
  {Swarttouw}(2010)}]{KOEKOEK:2010}%
  \BibitemOpen
  \bibfield  {author} {\bibinfo {author} {\bibfnamefont {R.}~\bibnamefont
  {Koekoek}}, \bibinfo {author} {\bibfnamefont {P.}~\bibnamefont {Lesky}}, \
  and\ \bibinfo {author} {\bibfnamefont {R.}~\bibnamefont {Swarttouw}},\
  }\href@noop {} {\emph {\bibinfo {title} {Hypergeometric Orthogonal
  Polynomials and Their q-Analogues}}},\ \bibinfo {edition} {1st}\ ed.,\ 185
  Springer Monographs in Mathematics\ (\bibinfo  {publisher}
  {Springer-Verlag},\ \bibinfo {address} {Berlin},\ \bibinfo {year}
  {2010})\BibitemShut {NoStop}%
\bibitem [{\citenamefont {Prudnikov}, \citenamefont {Brychkov},\ and\
  \citenamefont {Marichev}(1990)}]{PRUDNIKOV2}%
  \BibitemOpen
  \bibfield  {author} {\bibinfo {author} {\bibfnamefont {A.}~\bibnamefont
  {Prudnikov}}, \bibinfo {author} {\bibfnamefont {Y.}~\bibnamefont {Brychkov}},
  \ and\ \bibinfo {author} {\bibfnamefont {O.}~\bibnamefont {Marichev}},\
  }\href@noop {} {\emph {\bibinfo {title} {Integrals and Series. Special
  Functions}}},\ Vol.~\bibinfo {volume} {2}\ (\bibinfo  {publisher} {Gordon and
  Breach Sci. Publ.},\ \bibinfo {address} {New York},\ \bibinfo {year}
  {1990})\BibitemShut {NoStop}%
\bibitem [{\citenamefont {Magnus}, \citenamefont {Oberhettinger},\ and\
  \citenamefont {Soni}(1966)}]{MAGNUS}%
  \BibitemOpen
  \bibfield  {author} {\bibinfo {author} {\bibfnamefont {W.}~\bibnamefont
  {Magnus}}, \bibinfo {author} {\bibfnamefont {F.}~\bibnamefont
  {Oberhettinger}}, \ and\ \bibinfo {author} {\bibfnamefont {R.}~\bibnamefont
  {Soni}},\ }\href@noop {} {\emph {\bibinfo {title} {Formulas and Theorems for
  the Special Functions of Mathematical Physics}}}\ (\bibinfo  {publisher}
  {Springer-Verlag},\ \bibinfo {year} {1966})\BibitemShut {NoStop}%
\end{thebibliography}%

\end{document}